\documentclass[twoside]{article}

\usepackage{amsmath,amsthm,amssymb,latexsym}
\usepackage[v2,tips]{xy}

\newcommand{\ip}[2]{{\langle {#1} , {#2} \rangle}}

\newcommand{\proten}{{\widehat{\otimes}}}

\newcommand{\mc}[1]{\mathcal{#1}}
\newcommand{\inten}{{\check{\otimes}}}
\newcommand{\id}{\operatorname{id}}

\theoremstyle{plain}%
\newtheorem{proposition}{Proposition}[section]%
\newtheorem{theorem}[proposition]{Theorem}%
\newtheorem{corollary}[proposition]{Corollary}%
\newtheorem{lemma}[proposition]{Lemma}%
\theoremstyle{definition}%
\newtheorem{definition}[proposition]{Definition}%

\newcommand{\ap}{{\operatorname{AP}}}
\newcommand{\wap}{{\operatorname{WAP}}}
\newcommand{\vnten}{\overline\otimes}

\begin{document}

\large
\title{Weakly almost periodic functionals on the measure algebra}
\author{Matthew Daws}
\maketitle

\begin{abstract}
It is shown that the collection of weakly almost periodic functionals on the convolution algebra of a commutative Hopf von Neumann algebra is a C$^*$-algebra.  This implies that the weakly almost periodic functionals on $M(G)$, the measure algebra of a locally compact group $G$, is a C$^*$-subalgebra of $M(G)^* = C_0(G)^{**}$.  The proof builds upon a factorisation result, due to Young and Kaiser, for weakly compact module maps.  The main technique is to adapt some of the theory of corepresentations to the setting of general reflexive Banach spaces.


\end{abstract}

\section{Introduction}

For a topological (semi-)group $G$, the space of \emph{weakly almost periodic}
functions on $G$ is the subspace of $C(G)$ consisting of those $f\in C(G)$
such that the (left) translates of $f$ form a relatively weakly-compact subset
of $C(G)$.  We denote this set by $\wap(G)$.  Then $\wap(G)$ is a unital C$^*$-subalgebra
of $C(G)$, say with character space $G^\wap$.  By continuity, we can extend the
product from $G$ to $G^\wap$, turning $G^\wap$ into a compact semigroup whose product
is separately continuous, a \emph{semitopological semigroup}.  Indeed, $G^\wap$
is the universal semitopological semigroup compactification of $G$.  See \cite{BJM}
or \cite{HS} for further details.

Now suppose that $G$ is a locally compact group, so we may form the Banach space
$L^1(G)$, which becomes a Banach algebra with the convolution product.  Then $L^\infty(G)$,
as the dual of $L^1(G)$, naturally becomes an $L^1(G)$-bimodule.  We define the
space of \emph{weakly almost periodic functionals} on $L^1(G)$, denoted by $\wap(L^1(G))$,
to be the collection of $f\in L^1(G)$ such that the map
\[ R_f:L^1(G) \rightarrow L^\infty(G), \quad a\mapsto a\cdot f \qquad (a\in L^1(G)), \]
is weakly-compact.  This is equivalent to the map $L_f:L^1(G)\rightarrow L^\infty(G),
a\mapsto f\cdot a$ being weakly-compact.
\"Ulger showed in \cite{ulger} that $\wap(L^1(G)) = \wap(G)$, where $C(G)$ is
naturally identified with a subspace of $L^\infty(G)$.  This fact also follows easily
from \cite[Lemma~6.3]{wong}, using the fact that if a set is relatively weakly compact,
then the weak- and weak$^*$-topology closures coincide.  Both these papers use simple
bounded approximate identity arguments.

The definition of $\wap(L^1(G))$ obviously generalises to any Banach algebra $\mc A$.
In general, we can say little about $\wap(\mc A)$, except for some interesting
links with the Arens products, see \cite{daws} and references therein, or
\cite[Section~1.4]{palmer} and \cite[Theorem~2.6.15]{dales}.
However, motivated by
the above example, we might expect that when $\mc A$ has a large amount of structure,
$\wap(\mc A)$ also might have extra structure.  In this paper, we shall investigate
$\wap(M(G))$, where $M(G)$ is the measure algebra over a locally compact group $G$.
In particular, we shall show that $\wap(M(G))$ is a C$^*$-subalgebra of $M(G)^*$,
where $M(G)$ is identified as the dual of $C_0(G)$, so that $M(G)^* = C_0(G)^{**}$ is
a commutative von Neumann algebra.

The central idea is to develop a theory of corepresentations on reflexive Banach spaces,
for commutative Hopf von Neumann algebras.  Our theory exactly replicates that for
Hilbert spaces, but care needs to be taken to ensure that everything works in our more
general setting.

The connection between weakly almost periodic functionals and
representations of Banach algebras goes back to Young, \cite{young}, and Kaiser,
\cite{kai}.  For $L^1(G)$, there is a correspondence between (non-degenerate)
representations of $L^1(G)$ and representations of $G$.  Using Young and Kaiser's work,
it is easy to see that weakly almost periodic functionals on $L^1(G)$ correspond to
coefficient functionals for representations of $G$ on reflexive spaces.  Then
multiplication of functions in $L^\infty(G)$ corresponds to tensoring representations.
The existence of reflexive tensor products (see \cite{as} for example) hence shows
that the product of two weakly almost periodic functionals is again weakly almost
periodic.  Of course, for $L^1(G)$, it is far easier to use \"Ulger's result, and then
argue directly that $\wap(G)$ is an algebra (which follows from Grothendieck's criteria
for weak compactness, see \cite{BJM}).  For $M(G)$, while $M(G)=L^1(X)$ for some
measure space $X$, we do not have that $X$ is a (semi)group, and so we turn to
corepresentations, which work with the algebra $M(G)^*$ directly.

The structure of the paper is as follows.  We first introduce some notions from the theory
of tensor products of Banach spaces, in particular the projective and injective tensor
norms.  We then define what a (commutative) Hopf von Neumann algebra is, and show
carefully that $M(G)$ (as well as $L^1(G)$) fits into this abstract framework.  For the
rest of the paper, we work with commutative Hopf von Neumann algebras, the results for
$M(G)$ (and, indeed, $L^1(G)$) being immediate corollaries.  As an immediate application,
we make a quick study of almost periodic functionals.  We then turn our attention to
weakly almost periodic functionals, and build a theory of corepresentations on
reflexive Banach spaces.  The final application is then obtained by checking that the usual way
of tensoring corepresentations still works in this more general setting.

For an introduction to quantum groups from a functional analysis viewpoint, \cite{ku},
or the pair of articles \cite{kt1} and \cite{kt2}, are very readable.  A good starting
point for details about (weakly) almost periodic functionals on general Banach algebras is
\cite{dul}.

A few notes on notion.  We generally follow \cite{dales} for details about Banach
algebras.  We write $E^*$ for the dual of a Banach space $E$, and use the dual pairing
notation $\ip{\mu}{x} = \mu(x)$, for $\mu\in E^*$ and $x\in E$.
We write $\mc B(E,F)$ for the collection of bounded linear maps from $E$ to $F$, we write
$\mc B(E,E)=\mc B(E)$, and we write $T^*$ for the linear adjoint of an operator $T$.
%We write $\kappa_E:E\rightarrow E^{**}$ for the canonical map from a Banach space $E$
%to its bidual.

\smallskip\textbf{Acknowledgements:}
The author would like to thank Garth Dales, for bringing this problem
to his attention, and for careful proofreading.  Thanks to Tony Lau for providing
the reference \cite{wong}.

\section{Hopf von Neumann algebras}

We start by recalling some elementary definitions and facts from the theory
of tensor products of Banach spaces.  We refer the reader to the books \cite{ryan} and
\cite{df}, or \cite[Chapter~VIII]{du}, for further details.

Let $E$ and $F$ be Banach
spaces.  The \emph{projective tensor norm}, $\|\cdot\|_\pi$, on $E\otimes F$ is defined by
\[ \|\tau\|_\pi = \inf\Big\{ \sum_{k=1}^n \|x_k\| \|y_k\| : \tau=\sum_{k=1}^n x_k
\otimes y_k \Big\} \qquad (\tau\in E\otimes F). \]
Then $E\proten F$, the \emph{projective tensor product of $E$ and $F$}, is the
completion of $E\otimes F$ with respect to $\|\cdot\|_\pi$.  The projective tensor
product has the property that any bounded, bilinear map $\psi: E\times F\rightarrow G$
admits a unique bounded linear extension $\tilde\psi: E\proten F\rightarrow G$, with
$\|\tilde\psi\| = \|\psi\|$.  For measure spaces $X$ and
$Y$, we have that $L^1(X)\proten L^1(Y) = L^1(X\times Y)$.  We identify $(E\proten F)^*$
with $\mc B(E,F^*)$ under the dual pairing
\[ \ip{T}{x\otimes y} = \ip{T(x)}{y} \qquad (T\in\mc B(E,F^*), x\in E, y\in F), \]
and using linearity and continuity.

The \emph{injective tensor norm}, $\|\cdot\|_\epsilon$, on $E\otimes F$ is defined by
regarding $E\otimes F$ as a subspace of $\mc B(E^*,F)$, where $\tau = \sum_{k=1}^n
x_k \otimes y_k$ induces the finite-rank operator
\[ E^*\rightarrow F, \quad \mu\mapsto \sum_{k=1}^n \ip{\mu}{x_k} y_k. \]
Then $E\inten F$, the \emph{injective tensor product of $E$ and $F$}, is the
completion of $E\otimes F$ with respect to $\|\cdot\|_\epsilon$.  For locally
compact Hausdorff spaces $K$ and $L$, we have that $C_0(K) \inten C_0(L) =
C_0(K\times L)$.  We write $\mc A(E,F)$ for the closure of the finite-rank operators
from $E$ to $F$; these are the \emph{approximable operators from $E$ to $F$}.
Then, almost by definition, we have that $\mc A(E,F) = E^* \inten F$.

There is a canonical norm-decreasing map $E\proten F\rightarrow E\inten F$.
By taking the adjoint, we get an injective contraction $(E\inten F)^* \rightarrow
\mc B(E,F^*)$.  The image, equipped with the norm induced by $(E\inten F)^*$,
is the space of \emph{integral operators}, $\mc I(E,F^*)$.
The map $E^*\proten E\rightarrow E^*\inten E$ is injective if and only if $E$ has the
\emph{approximation property}.  We can regard $E^*\proten E$ as a subspace of
$\mc I(E^*) = (E^*\inten E)^* = \mc A(E)^*$ by
\[ \ip{\mu\otimes x}{T} = \ip{\mu}{T(x)} \qquad (T\in\mc A(E), \mu\otimes x\in
E^*\proten E). \]
Similarly, we can regard $E\proten F$ as a subspace of $\mc I(E^*,F)$; here we use
that fact that $\mc I(E^*,F)$ is isometrically a subspace of $\mc I(E^*,F^{**})
= (E^*\inten F^*)^* = \mc A(E,F^*)^*$.  We say that $E$ has the \emph{metric approximation
property} if and only if the map $E^*\proten E\rightarrow\mc I(E^*)$ is an isometry
onto its range, or equivalent, $E\proten F \rightarrow \mc I(E^*,F)$ is an
isometry onto its range, for all $F$.  There are characterisations of the (metric)
approximation property in terms of finite-rank approximations of the identity on compact
sets.  We have that $C_0(K)$ and $L^1(X)$ have the metric approximation property for
all $K$ and $X$.

\subsection{Commutative Hopf von Neumann algebras}

A \emph{Hopf von Neumann algebra} is a von Neumann algebra $\mc M$ equipped with
a \emph{coproduct} $\Delta:\mc M\rightarrow \mc M\vnten\mc M$.  Here $\vnten$ denotes
the von Neumann tensor product.  This means that
$\Delta$ is a normal $*$-homomorphism, and that $(\Delta\otimes\id)\Delta =
(\id\otimes\Delta)\Delta$, that is, $\Delta$ is coassociative.

We shall concentrate on the case where $\mc M$ is commutative, so that $\mc M = L^\infty(X)$
for some measure space $X$.  Then $\mc M \vnten \mc M = L^\infty(X\times X)$, and so, as
$\Delta$ is normal, it drops to give a contractive map $\Delta_*:L^1(X\times X) =
L^1(X) \proten L^1(X) \rightarrow L^1(X)$.  Hence $\Delta_*$ induces a contractive bilinear
map $L^1(X)\times L^1(X) \rightarrow L^1(X)$.  Then $\Delta$ being coassociative is
equivalent to $\Delta_*$ being associative.

As both $\mc M$ and $\mc M_*$ are Banach algebras, we have natural module actions
of $\mc M$ on $\mc M_*$ and of $\mc M_*$ on $\mc M$.  For the action of $\mc M$ on
$\mc M_*$, we shall, for example, write $F\cdot a\in\mc M_*$ for $F\in\mc M$ and
$a\in\mc M_*$.  For the action of $\mc M_*$ on $\mc M$, we shall always explicitly
invoke the map $\Delta_*$ or $\Delta$.

For an example of a commutative Hopf von Neumann algebra, let $G$ be a locally compact
group, and consider the algebra $L^\infty(G)$
equipped with the coproduct $\Delta$ defined by
\[ \Delta(f)(s,t) = f(st) \qquad (f\in L^\infty(G), s,t\in G). \]
Then $\Delta_*$ induces the usual convolution product on $L^1(G)$.

A slightly less well-known example is furnished by $M(G)$.  As $M(G) = C_0(G)^*$, we see that
$M(G)$ is the predual of the commutative von Neumann algebra $C_0(G)^{**}$.  As such,
$M(G)^* = L^\infty(X)$ for some measure space $X$ (see \cite[Chapter~III]{tak}), and
so by the uniqueness of preduals, $M(G) = L^1(X)$.  Let $\Phi$ be the canonical coproduct
on $C_0(G)$, so that $\Phi$ is the $*$-homomorphism $C_0(G)\rightarrow C(G\times G)$ defined
by \[ \Phi(f)(s,t) = f(st) \qquad (f\in C_0(G), s,t\in G). \]
We identify $C(G\times G)$, the space of bounded continuous functions on $G\times G$,
with the \emph{multiplier algebra} of $C_0(G\times G)$,
and hence (see \cite[Chapter~III, Section~6]{tak}) we may identify $C(G\times G)$ with
\[ \{ x\in C_0(G\times G)^{**} : fx,xf\in C_0(G\times G) \ (f\in C_0(G\times G)) \}. \]
We can hence regard $\Phi$ as a $*$-homomorphism $C_0(G) \rightarrow C_0(G\times G)^{**}$.
%Hence there is a unique extension $\tilde\Phi: C_0(G)^{**} \rightarrow
%C_0(G\times G)^{**}$, compare \cite[Chapter~III, Lemma~2.2]{tak}.

We claim that $M(G) \proten M(G)$ is, isometrically, a subspace of $M(G\times G)
= C_0(G\times G)^*$.  From the above, we can identify $M(G\times G)$ with
$\mc I(C_0(G), M(G))$.  As $M(G)$ has the metric approximation property, we see
that $M(G) \proten M(G)$ is isometrically a subspace of $\mc I(M(G)^*,M(G))$, or
equivalently, by properties of the integral operators, isometrically a subspace of
$\mc I(C_0(G),M(G))$, as required.

Alternatively, for any C$^*$-algebra $\mc A$,
we could define a norm on $\mc A^*\otimes\mc A^*$ by embedding $\mc A^*\otimes\mc A^*$
into $(\mc A\otimes_{\min}\mc A)^*$.  This induces the \emph{operator space projective
tensor norm}, see \cite[Chapter~7]{er}, and as $\mc A$ has the \emph{minimal} operator
space structure, it follows that $\mc A^*$ has the \emph{maximal} structure, and
so this norm agrees with the (Banach space) projective tensor norm.

Hence $L^\infty(X)\vnten L^\infty(X) = (M(G)\proten M(G))^*$ is a quotient of
$M(G\times G)^* = C_0(G\times G)^{**}$.  We claim that this quotient map is a 
$*$-homomorphism, for which it suffices to check that the kernel
\[ \{ \tau\in C_0(G\times G)^{**} : \ip{\tau}{\mu\otimes\lambda}=0 \ (\mu,\lambda
\in M(G)) \} \]
is an ideal.  Let $\mu,\lambda\in M(G)$, let $g\in C_0(G\times G)$, and let
$f=f_1\otimes f_2 \in C_0(G) \otimes C_0(G)$.  Then
\[ \ip{(\mu\otimes\lambda)\cdot f}{g} = \int f_1(s) f_2(t) g(s,t) \ d\mu(s) \ d\lambda(t), \]
so $(\mu\otimes\lambda)\cdot f$ is the measure $\mu\cdot f_1 \otimes \lambda\cdot f_2
\in M(G) \otimes M(G)$.  By continuity, we see that $(\mu\otimes\lambda)\cdot f
\in M(G) \proten M(G)$ for any $f\in C_0(G\times G)=C_0(G)\inten C_0(G)$.
Let $\tau$ be in the kernel, so that
\[ \ip{\tau\cdot(\mu\otimes\lambda)}{f} = \ip{\tau}{(\mu\otimes\lambda)\cdot f} = 0, \]
as $\tau$ kills $M(G) \proten M(G)$.  Thus $\tau\cdot(\mu\otimes\lambda)=0$ in
$C_0(G\times G)^*$.  So, for $\sigma\in C_0(G\times G)^{**}$, we see that
\[ \ip{\sigma\tau}{\mu\otimes\lambda} = \ip{\sigma}{\tau\cdot(\mu\otimes\lambda)} = 0, \]
so that $\sigma\tau$ lies in the kernel.

Hence we have the following chain of $*$-homomorphisms
\[ \xymatrix{ C_0(G) \ar[r]^{\Phi} & C_0(G\times G)^{**}
\ar[r] & L^\infty(X\times X) = C_0(G)^{**} \vnten C_0(G)^{**}, } \]
say, giving rise to a $*$-homomorphism $\Delta_0:C_0(G)\rightarrow L^\infty(X\times X)$.
There is hence a canonical extension (compare \cite[Chapter~III, Lemma~2.2]{tak})
$\Delta:L^\infty(X) \rightarrow L^\infty(X\times X)$, which is a normal $*$-homomorphism.
Indeed, the preadjoint $\Delta_*:M(G)\proten M(G)\rightarrow M(G)$ is defined by the
chain of maps
\[ \xymatrix{ M(G)\proten M(G) = L^1(X\times X) \ar[r] & L^\infty(X\times X)^*
\ar[r]^{\Delta_0^*} & C_0(G)^* = M(G). } \]
Then, for $\mu,\lambda\in M(G)$ and $f\in C_0(G)$, we see that
\begin{align*} \ip{\Delta_*(\mu\otimes\lambda)}{f} &=
\ip{\Delta_0(f)}{\mu\otimes\lambda} = \int_{G\times G} f(st) \ d\mu(s) \ d\lambda(t),
\end{align*}
so $\Delta_*$ induces the usual convolution product on $M(G)$.  We have hence shown
that $M(G)^*$ is a commutative Hopf von Neumann algebra.
Notice that throughout, we have actually only used the fact that $G$ is a
locally compact semigroup.

For a recent survey on measure algebras, see \cite{DLS}, where the authors view
$M(G)$ as a \emph{Lau algebra} (see \cite{Lau}).

\subsection{Almost periodic functionals}

For a Banach algebra $\mc A$, a functional $\mu\in\mc A^*$ is \emph{almost periodic}
if the map
\[ R_\mu: \mc A\rightarrow\mc A^*, \quad a\mapsto a\cdot\mu \qquad (a\in\mc A), \]
is compact.  We denote the collection of almost periodic functionals by $\ap(\mc A)$.
Then it is easy to see that $\ap(\mc A)$ is a closed subspace of $\mc A^*$.  Using
the viewpoint of Hopf von Neumann algebras, it is easy to see that
$\ap(M(G))$ is a C$^*$-algebra.

\begin{theorem}\label{ap_case}
Let $(L^\infty(X),\Delta)$ be a commutative Hopf von Neumann algebra, so that $L^1(X)$
becomes a Banach algebra.  Then $\ap(L^1(X))$ is a C$^*$-subalgebra of $L^\infty(X)$.
\end{theorem}
\begin{proof}
Let $F\in\ap(L^1(X))\subseteq L^\infty(X)$.  For $f\in L^1(X)$, we shall write
$f^*$ for the pointwise complex-conjugation of $f$, so that $f\mapsto f^*$ is
the preadjoint of the involution on $L^\infty(X)$.  We see that for $f,g\in L^1(X)$,
\begin{align*} \ip{R_{F^*}(f)}{g} &= \ip{F^*}{\Delta_*(g\otimes f)}
= \ip{\Delta(F)^*}{g\otimes f} = \ip{F}{\Delta_*(g^*\otimes f^*)} \\
&= \ip{R_F(f^*)}{g^*} = \ip{R_F(f^*)^*}{g}, \end{align*}
so we conclude that $R_{F^*}$ is compact if and only if $R_F$ is compact.  Hence
$\ap(L^1(X))$ is $*$-closed.

We claim that $R_F = \Delta(F)^*\kappa_{L^1(X)}$.  Indeed, for $f,g\in L^1(X)$, we have that
\begin{align*} \ip{R_F}{f\otimes g} &= \ip{R_F(f)}{g} = \ip{F}{\Delta_*(g\otimes f)}
= \ip{\Delta(F)}{g\otimes f} \\ &= \ip{\Delta(F)(g)}{f}
= \ip{\Delta(F)\kappa_{L^1(X)}(f)}{g}. \end{align*}
So $\Delta(F) = R_F^*\kappa_{L^1(X)}$, and hence $R_F$ is compact if and only if
$\Delta(F)$ is compact.
As $L^\infty(X)$ has the approximation property, it follows that $\Delta(F)$ is compact
if and only if \[ \Delta(F)\in\mc A(L^1(X),L^\infty(X)) = L^\infty(X) \inten L^\infty(X)
= L^\infty(X) \otimes_{\min} L^\infty(X) \subseteq L^\infty(X) \vnten L^\infty(X). \]
Thus, if $F,G\in\ap(L^1(X))$, then $\Delta(F),\Delta(G) \in
L^\infty(X)\otimes_{\min} L^\infty(X)$, and so $\Delta(FG) = \Delta(F)\Delta(G)
\in L^\infty(X)\otimes_{\min} L^\infty(X)$, as $L^\infty(X)\otimes_{\min} L^\infty(X)$
is an algebra.  Hence $FG\in\ap(L^1(X))$, as required.
\end{proof}

\begin{corollary}
For a locally compact group $G$, $\ap(M(G))$ is a C$^*$-subalgebra of $M(G)^*$.
\end{corollary}

\section{Weakly almost periodic functionals}

We shall make use of vector valued $L^p$ spaces; for a measure space $X$, a Banach space
$E$, and $1\leq p<\infty$, we write $L^p(X,E)$ for the space of (classes of almost
everywhere equal) Bochner $p$-integrable functions from $X$ to $E$.  Then $L^p(X)
\otimes E$ naturally maps into $L^p(X,E)$ with dense range, inducing a norm
$\Delta_p$ on $L^p(X) \otimes E$.  This norm is studied in \cite[Chapter~7]{df}.
We have that $L^1(X) \proten E = L^1(X,E)$, so that $\Delta_1 = \|\cdot\|_\pi$, the
projective tensor norm.

It is worth noting that $\Delta_p$ is not a \emph{tensor norm}, as $T\in\mc B(L^p(X))$
may fail to extend to a bounded map $T\otimes\id: L^p(X,E)\rightarrow L^p(X,E)$.
%We shall return to this idea later.
However, note that for $F\in L^\infty(X)$,
then denoting also by $F$ the multiplication operator on $L^p(X)$, it is elementary
that $F\otimes\id$ is bounded, with norm $\|F\|$, on $L^p(X,E)$.  The norm
$\Delta_p$ does satisfy the estimates
\[ \|\tau\|_\epsilon \leq \Delta_p(\tau) \leq \|\tau\|_\pi
\quad (\tau \in L^p(X)\otimes E), \]
so in particular, $\Delta_p(f\otimes x) = \|f\| \|x\|$ for $f\in L^p(X)$ and $x\in E$.

We shall henceforth restrict to the case where $E$ is reflexive.  Then $E^*$ has the
Radon-Nikod\'ym property, and so $L^p(X,E)^* = L^{p'}(X,E^*)$ for $1<p<\infty$,
where $1/p' = 1-1/p$, see \cite[Appendix~D]{df}, or \cite{du}, for further details.
We stress that even when $p=2$, the dual pairing between $L^2(X,E)$ and $L^2(X,E^*)$
is always \emph{bilinear} and not sesquilinear.

\begin{lemma}
Let $E$ be a reflexive Banach space, and let $X$ be a measure space.  The map
\[ \Lambda:\big( L^2(X) \otimes E^* \big) \times \big( L^2(X) \otimes E \big)
\rightarrow L^1(X) \otimes \big(E^*\proten E\big);
\big( f\otimes\mu , g\otimes x\big) \mapsto fg \otimes (\mu\otimes x) \]
extends to a metric surjection
\[ \Lambda:L^2(X,E^*) \proten L^2(X,E) \rightarrow L^1\big(X,E^*\proten E\big) =
L^1(X) \proten E^* \proten E. \]
Here $fg$ denotes the pointwise product, so the Cauchy-Schwarz inequality shows that
$fg\in L^1(X)$ for $f,g\in L^2(X)$.
\end{lemma}
\begin{proof}
Let $F\in L^2(X,E^*)$ and $G\in L^2(X,E)$ be simple functions, so that
there exists a disjoint partition of $X$, say $(X_k)_{k=1}^n$, and $(x_k)_{k=1}^n
\subseteq E$ and $(\mu_k)_{k=1}^n \subseteq E^*$ with
\[ F = \sum_{k=1}^n \chi_{X_k} \otimes \mu_k, \quad
G = \sum_{k=1}^n \chi_{X_k} \otimes x_k. \]
Here we write $\chi_{X_k}$ for the indicator function of $X_k$.
Hence we see that
\[ \Lambda(F\otimes G) = \sum_{k=1}^n \chi_{X_k} \otimes (\mu_k\otimes x_k), \]
which has norm
\begin{align*}
\sum_{k=1}^n |X_k| \|\mu_k\otimes x_k\| &\leq
\Big( \sum_{k=1}^n |X_k| \|\mu_k\|^2 \Big)^{1/2}
\Big( \sum_{k=1}^n |X_k| \|x_k\|^2 \Big)^{1/2} = \|F\| \|G\|.
\end{align*}
As the simple functions are dense in $L^2(X,E)$, respectively, $L^2(X,E^*)$, we conclude
that the map $\Lambda:L^2(X,E^*) \times L^2(X,E) \rightarrow L^1(X,E^*\proten E)$
is a contraction, and so extends to a contraction $L^2(X,E^*) \proten L^2(X,E)
\rightarrow L^1(X,E^*\proten E)$.

As $E$ is reflexive, we may identify $(E^*\proten E)^*$ with $\mc B(E)$.
Hence $\Lambda^*$ is a map $\mc B(L^1(X),\mc B(E)) \rightarrow \mc B(L^2(X,E))$, say
$\pi \mapsto W$, where
\[ \ip{f\otimes\mu}{W(g\otimes x)} = \ip{\mu}{\pi(fg)(x)}
\qquad (f,g\in L^2(X), \mu\in E^*, x\in E). \]
By a suitable choice of $f,g,x$ and $\mu$, we see that $\|W\| \geq \|\pi\|$,
and so we conclude that actually $\|W\| = \|\pi\|$.  Hence $\Lambda^*$
is an isometry, so $\Lambda$ must be a metric surjection, as required.
\end{proof}

For $F\in L^\infty(G)$ and $T\in\mc B(E)$, we see that $F\otimes T$ extends to a bounded
linear map on $L^2(X,E)$.  Let $L^\infty(X) \vnten \mc B(E)$ be the weak$^*$-closure
of $L^\infty(X)\otimes \mc B(E)$ inside $\mc B(L^2(X,E))$.  This is then a \emph{dual
Banach algebra}, that is, multiplication in $L^\infty(X) \vnten \mc B(E)$ is
separately weak$^*$-continuous.  See \cite[Section~8]{daws}, where similar ideas
are explored.

\begin{proposition}\label{image_of_lambda_ad}
The above lemma isometrically identifies $\mc B(L^1(X),\mc B(E))$ with a subspace
of $\mc B(L^2(X,E))$, under the mapping $\Lambda^*$.  The image of $\Lambda^*$
is precisely $L^\infty(X) \vnten \mc B(E)$.
\end{proposition}
\begin{proof}
Standard Banach space theory shows that the image of $\Lambda^*$ is equal to
\[ (\ker\Lambda)^\perp = \big\{ T\in\mc B(L^2(X,E)) : \ip{T}{\tau}=0 \
(\tau\in L^2(X,E^*) \proten L^2(X,E), \Lambda(\tau)=0) \big\}. \]
Hence the image of $\Lambda^*$ is weak$^*$-closed.  Notice that $L^\infty(X) \vnten
\mc B(E)$ is equal to $Z^\perp$, where
\[ Z = \big\{ \tau\in L^2(X,E^*) \proten L^2(X,E) : \ip{F\otimes S}{\tau}=0
\ (F\in L^\infty(X), S\in\mc B(E)) \big\}. \]
Hence we need to show that $\ker\Lambda = Z$.

Let $F\in L^\infty(X)$ and $S\in\mc B(E)$.  Then let $T=F\otimes S\in \mc B(L^2(X,E))$,
and let $\pi:L^1(X) \rightarrow \mc B(E)$ be the rank-one operator induced by
$F\otimes S$, that is, $\pi(a)=\ip{F}{a}S$ for $a\in L^1(X)$.
Then $\Lambda^*(\pi) = T$, from which it follows that $\ker\Lambda\subseteq Z$.

As $L^1(X)$ has the approximation property, for each non-zero $\sigma \in L^1(X)
\proten (E^*\proten E)$, there exists $F\in L^\infty(X)$ and $S\in\mc B(E)$ with
$\ip{F\otimes S}{\sigma}\not=0$.  Hence, if $\tau\in L^2(X,E^*) \proten L^2(X,E)$
is such that $\sigma = \Lambda(\tau) \not=0$, then there exists $T\in L^\infty(X)
\otimes \mc B(E)$ with $0\not=\ip{T}{\sigma}=\ip{\Lambda^*(T)}{\tau}$.  This
shows that $Z \subseteq \ker\Lambda$.
\end{proof}

Informally, the above proposition allows us to write
\[ \mc B(L^1(X),\mc B(E)) = \big( L^1(X) \proten (E^*\proten E) \big)^*
= L^\infty(X) \vnten \mc B(E), \]
which is reminiscent of the operator space projective tensor result that
$(\mc M_* \proten \mc N_*)^* = \mc M \vnten \mc N$, for von Neumann algebras $\mc M$
and $\mc N$, see \cite[Theorem~7.2.4]{er}.  The important point for us is that we
have turned $\mc B(L^1(X),\mc B(E))$ into an algebra.  It is multiplication in this
algebra which will ultimately give rise to the multiplication of weakly almost periodic
functionals in $L^\infty(X)$.

Now let $L^\infty(X)$ be a Hopf von Neumann algebra, so it admits a coproduct $\Delta$.
We have a map
\[ \Delta_*\otimes\id: L^1(X\times X) \proten (E^*\proten E) \rightarrow
L^1(X) \proten (E^*\proten E), \]
whose adjoint, which we denote by $\Delta\otimes\id$, is a map
\[ \Delta\otimes\id: L^\infty(X) \vnten \mc B(E) \rightarrow
L^\infty(X\times X) \vnten \mc B(E), \]
where, of course, $L^\infty(X\times X) \vnten \mc B(E)$ is a subalgebra of
$\mc B(L^2(X\times X,E))$.

\begin{lemma}\label{delta_homo}
With notation as above, $(\Delta\otimes\id)$ is a homomorphism.
%That is, for $U,V\in L^\infty(X) \vnten \mc B(E)$, we have that
%$(\Delta\otimes\id)(UV) = ( (\Delta\otimes\id)U )( (\Delta\otimes\id)V )$.
\end{lemma}
\begin{proof}
Let $F\in L^\infty(X)$ and $a,b\in L^1(X)$.  As $\Delta$ is a homomorphism,
it follows that
\[ \Delta_*\big( \Delta(F) \cdot (a\otimes b) \big) =
F \cdot \Delta_*(a\otimes b). \]

Let $U\in L^\infty(X)\vnten\mc B(E)$, let $T\in\mc B(E)$, and let $V=F\otimes T$.
For $\tau\in E^*\proten E$, we see that
\begin{align*} \ip{ (\Delta\otimes\id)(UV) }{a\otimes b\otimes\tau}
&= \ip{U(F\otimes T)}{\Delta_*(a\otimes b)\otimes\tau}
= \ip{U}{F\cdot\Delta_*(a\otimes b)\otimes T\cdot\tau} \\
&= \ip{U}{ \Delta_*( \Delta(F) \cdot (a\otimes b)) \otimes T\cdot\tau} \\
&= \ip{(\Delta\otimes\id)U}{\Delta(F) \cdot (a\otimes b) \otimes T\cdot\tau} \\
&= \ip{((\Delta\otimes\id)U)((\Delta\otimes\id)V)}{a\otimes b \otimes\tau}. \end{align*}
By linearity, we conclude that $(\Delta\otimes\id)(UV) =
( (\Delta\otimes\id)U )( (\Delta\otimes\id)V )$ for all $U\in L^\infty(X)\vnten\mc B(E)$
and $V\in L^\infty(X) \otimes \mc B(E)$.  By weak$^*$-continuity, this must also hold
for $V\in L^\infty(X) \vnten \mc B(E)$.
\end{proof}

We now wish to adapt \emph{leg numbering notation} to our setup.
Given $W\in\mc B(L^2(X,E))$, define $W_{23} \in \mc B(L^2(X\times X,E))$ by
\[ W_{23}(f_1\otimes f_2\otimes x) = f_1 \otimes W(f_2\otimes x)
\qquad (f_1,f_2\in L^2(X), x\in E). \]
Using the fact that $L^2(X\times X,E) = L^2(X, L^2(X,E))$, it is easy to see that
$W\mapsto W_{23}$ is a weak$^*$-continuous, isometric mapping.  If $W=
F\otimes S$ for some $F\in L^\infty(X)$ and $S\in\mc B(E)$, then clearly
$W_{23} = 1 \otimes F \otimes S \in L^\infty(X\times X) \otimes \mc B(E)$.  By
weak$^*$-continuity, we conclude that if $W\in L^\infty(X) \vnten \mc B(E)$, then
$W_{23} \in L^\infty(X\times X) \vnten \mc B(E)$.

Let $\chi:L^2(X\times X) \rightarrow L^2(X\times X)$ be the ``swap map'', defined on
elementary tensors by $\chi(f\otimes g) = g\otimes f$.  For
$W\in L^\infty(X\times X)\vnten \mc B(E)$, it is clear that $(\chi\otimes\id)W$
and $W(\chi\otimes\id)$ both also lie in $L^\infty(X\times X)\vnten \mc B(E)$.
For $W\in L^\infty(X)\vnten\mc B(E)$, we define $W_{13} = (\chi\otimes\id) W_{23}
(\chi\otimes\id) \in L^\infty(X\times X) \vnten \mc B(E)$.

\begin{theorem}\label{when_rep}
Let $(L^\infty(X),\Delta)$ be a Hopf von Neumann algebra, and let $E$ be a reflexive
Banach space.
Let $\pi:L^1(X)\rightarrow\mc B(E)$ be a bounded linear map, giving rise to
$W \in L^\infty(X) \vnten \mc B(E)$.  Then $\pi$ is a homomorphism, with respect
to $\Delta_*$, if and only if $(\Delta\otimes\id)W = W_{13} W_{23}$.
\end{theorem}
\begin{proof}
Let $f_1,f_2,g_1,g_2\in L^2(X)$, $\mu\in E^*$ and $x\in E$.  Then
\begin{align*} \ip{f_1\otimes f_2\otimes\mu}{(\Delta\otimes\id)W (g_1\otimes g_2\otimes x)}
&= \ip{\pi}{\Delta_*(f_1g_1 \otimes f_2g_2) \otimes (\mu\otimes x)} \\
&= \ip{\mu}{ \pi( \Delta_*(f_1g_1\otimes f_2g_2)) (x) }. \end{align*}

We now come to a proof where ``Sweedler notation'' would help greatly, but we should
perhaps, at least once, give a formal proof.  Informally, we shall ``pretend'' that
$W(g_2\otimes x) = h\otimes y$.  Then
\begin{align*} \ip{\mu}{\pi(f_1g_1)\pi(f_2g_2)(x)}
&= \ip{\pi(f_1g_1)^*(\mu)}{\pi(f_2g_2)(x)}
= \ip{f_2 \otimes \pi(f_1g_1)^*(\mu)}{W(g_2\otimes x)} \\
&= \ip{f_2}{h} \ip{\mu}{\pi(f_1g_1)(y)}
= \ip{f_2}{h} \ip{f_1\otimes\mu}{W(g_1\otimes y)} \\
&= \ip{f_1\otimes f_2\otimes\mu}{W_{13}(g_1\otimes h\otimes y)} \\
&= \ip{f_1\otimes f_2\otimes\mu}{W_{13}W_{23}(g_1\otimes g_2\otimes x)}, \end{align*}
which completes the proof.

To make this rigorous, for $\epsilon>0$, we can find a finite sum of
elementary tensors
$\sum_k h_k \otimes y_k \in L^2(X)\otimes E$ with $\|W(g_2\otimes x) -
\sum_k h_k\otimes y_k \| < \epsilon$.  Then
\[ \Big\| \ip{f_2 \otimes \pi(f_1g_1)^*(\mu)}{W(g_2\otimes x)} -
\sum_k \ip{f_2}{h_k} \ip{\mu}{\pi(f_1g_1)(y_k)} \Big\| < \epsilon\|f_2\|
\|\pi\| \|f_1\| \|g_1\| \|\mu\|, \]
and, as above,
\[ \sum_k \ip{f_2}{h_k} \ip{\mu}{\pi(f_1g_1)(y_k)} =
\sum_k \ip{f_1\otimes f_2\otimes\mu}{W_{13}(g_1\otimes h_k\otimes y_k)}, \]
so approximating again,
\begin{align*} \Big\|\sum_k & \ip{f_1\otimes f_2\otimes\mu}{W_{13}(g_1\otimes h_k\otimes y_k)}
- \ip{f_1\otimes f_2\otimes\mu}{W_{13}W_{23}(g_1\otimes g_2\otimes x)} \Big\|
\\ &< \epsilon \|W\|\|f_1\| \|f_2\| \|\mu\| \|g_1\|, \end{align*}
and so
\begin{align*} \big| \ip{\mu}{\pi(f_1g_1)\pi(f_2g_2)(x)} &  -
\ip{f_1\otimes f_2\otimes\mu}{W_{13}W_{23}(g_1\otimes g_2\otimes x)} \big| \\
&< 2\epsilon \|W\|\|f_1\| \|f_2\| \|\mu\| \|g_1\|. \end{align*}
As $\epsilon>0$ was arbitrary, the proof is complete.
\end{proof}

\subsection{Application to weakly almost periodic elements}

The following result was first shown by Young in \cite{young}, building upon \cite{DFJP},
and was recast in terms of the real interpolation method by Kaiser in \cite{kai}
(see also the similar arguments in \cite{daws}).

\begin{theorem}\label{wap_rep}
Let $\mc A$ be a Banach algebra, and let $\mu\in\mc A^*$.  The following are equivalent:
\begin{enumerate}
\item $\mu\in\wap(\mc A)$;
\item there exists a reflexive Banach space $E$, a contractive homomorphism $\pi:
\mc A\rightarrow\mc B(E)$, and $x\in E$, $\lambda\in E^*$ such that
\[ \ip{\mu}{a} = \ip{\lambda}{\pi(a)(x)} \qquad (a\in\mc A). \]
\end{enumerate}
\end{theorem}

We shall need a way of tensoring reflexive Banach spaces in a way that gives a reflexive
Banach space.  As we do not wish to get bogged down in the details of any one specific
way to do this, so we shall make an ad hoc definition.

\begin{definition}
Let $E$ and $F$ be reflexive Banach spaces, and let $\alpha$ be some norm on
$E\otimes F$ such that:
\begin{enumerate}
\item we have that $\|\tau\|_\epsilon \leq \alpha(\tau) \leq \|\tau\|_\pi$ for
  each $\tau\in E\otimes F$;
\item the completion of $E\otimes F$ with respect to $\alpha$ is a reflexive
  Banach space, say $E\proten_\alpha F$;
\item given $T\in\mc B(E)$ and $S\in\mc B(F)$, the map $T\otimes S:E\otimes F
  \rightarrow E\otimes F$ extends to a bounded operator on $E\proten_\alpha F$
  with norm $\|T\| \|S\|$.
\end{enumerate}
Then we say that $\|\cdot\|$ is a \emph{reflexive tensor norm} on $E\otimes F$.
\end{definition}

The existence of reflexive tensor norms is shown in \cite{as}, for example.

Let $E$ and $F$ be reflexive Banach spaces, and let $\alpha$ be a reflexive tensor norm
on $E\otimes F$.  As the map $E\proten F\rightarrow E \proten_\alpha F$ is
contractive with dense range, the adjoint $(E\proten_\alpha F)^* \rightarrow
\mc B(E,F^*)$ is injective.  We write $\mc B_{\alpha'}(E,F^*)$ for the image,
and equip it with the norm coming from $(E\proten_\alpha F)^*$, so we may write
$(E\proten_\alpha F)^* = \mc B_{\alpha'}(E,F^*)$.  For some norms $\alpha$, there
exists a \emph{dual norm} $\alpha'$, which is a reflexive tensor norm on
$E^* \otimes F^*$, such that $\mc B_{\alpha'}(E,F^*) = E^* \proten_{\alpha'} F^*$.
We shall, however, not have to assume this extra condition.  For us, it suffices
to note that as $\alpha$ dominates $\|\cdot\|_\epsilon$, there is a natural
embedding of $E^*\otimes F^*$ into $(E\proten_\alpha F)^*$.

\begin{theorem}
Let $(L^\infty(X),\Delta)$ be a commutative Hopf von Neumann algebra, and use $\Delta_*$
to turn $L^1(X)$ into a Banach algebra.  Then $\wap(L^1(X))$ is a C$^*$-subalgebra of
$L^\infty(X)$.
\end{theorem}
\begin{proof}
We know that $\wap(L^1(X))$ is a closed subspace of $L^\infty(X)$.  Exactly the same
argument as in the proof of Theorem~\ref{ap_case} shows that $\wap(L^1(X))$ is $*$-closed,
so it remains only to show that $\wap(L^1(X))$ is closed under multiplication.

Let $F_1,F_2\in\wap(L^1(X))$.  By Theorem~\ref{wap_rep}, for $i=1,2$, there exists a
reflexive Banach space $E_i$, a contractive homomorphism
$\pi_i:L^1(X)\rightarrow\mc B(E_i)$, $x_i\in E_i$ and $\mu_i\in E_i^*$ such that
\[ \ip{F_i}{a} = \ip{\mu_i}{\pi_i(a)(x_i)} \qquad (a\in L^1(X)). \]
Let $\alpha$ be a reflexive tensor norm on $E_1\otimes E_2$, and define
$\hat\pi_i:L^1(X) \rightarrow \mc B(E_1\proten_\alpha E_2)$, for $i=1,2$, by
\[ \hat\pi_1(a) = \pi_1(a)\otimes\id, \quad
\hat\pi_2(a) = \id\otimes\pi_2(a) \qquad (a\in L^1(X)). \]
Then $\hat\pi_1$ and $\hat\pi_2$ are contractive homomorphisms from $L^1(X)$ to
$\mc B(E_1\proten_\alpha E_2)$, and hence give rise, respectively, to $U,V\in L^\infty(X)
\vnten \mc B(E_1\proten_\alpha E_2)$, such that
\[ (\Delta\otimes\id)U = U_{13}U_{23}, \quad
(\Delta\otimes\id)V = V_{13}V_{23}. \]

Let $W = UV \in L^\infty(X) \vnten \mc B(E_1\proten_\alpha E_2)$, and let
$\pi:L^1(X) \rightarrow \mc B(E_1\proten_\alpha E_2)$ be induced by $W$.
By Lemma~\ref{delta_homo}, we see that
\[ (\Delta\otimes\id)W = ((\Delta\otimes\id)U)((\Delta\otimes\id)V)
= U_{13} U_{23} V_{13} V_{23}. \]
We also see that
\[ W_{13} W_{23} = U_{13} V_{13} U_{23} V_{23}. \]
We claim that $U_{23} V_{13} = V_{13} U_{23}$, from which it follows, from
Theorem~\ref{when_rep}, then $\pi$ is a homomorphism.

To prove the claim, we shall again deploy Sweedler notation: the argument can be
made rigorous in the same way as in the proof of Theorem~\ref{when_rep}.
Let $f_1,f_2,g_1,g_2\in L^2(X)$, $w_1\in E_1$, $w_2\in E_2$ and $T\in\mc B_{\alpha'}(E_1,
E_2^*) = (E_1\proten_\alpha E_2)^*$.  Suppose that $U^*(f_2\otimes T) = h_1\otimes S_1$,
so that for $k\in L^2(X)$, $z_1\in E_1$ and $z_2\in E_2$, we have that
\begin{align*} \ip{h_1\otimes S_1}{k\otimes z_1\otimes z_2}
&= \ip{f_2\otimes T}{U(k\otimes z_1\otimes z_2)}
= \ip{T}{\pi_1(f_2k)(z_1)\otimes z_2} \\&= \ip{T\pi_1(f_2k)(z_1)}{z_2}. \end{align*}
Similarly, suppose that $V^*(f_1\otimes T) = h_2\otimes S_2$, so that
\begin{align*} \ip{h_2\otimes S_2}{k\otimes z_1\otimes z_2}
&= \ip{f_1\otimes T}{V(k\otimes z_1\otimes z_2)}
= \ip{T}{z_1\otimes \pi_2(f_1k)(z_2)} \\&= \ip{T(z_1)}{\pi_2(f_1k)(z_2)}. \end{align*}
Thus we see that
\begin{align*}
&\ip{f_1\otimes f_2\otimes T}{U_{23} V_{13}(g_1\otimes g_2\otimes w_1\otimes w_2)} \\
&= \ip{f_1\otimes U^*(f_2\otimes T)}{(\chi\otimes\id\otimes\id)
   (g_2\otimes V(g_1\otimes w_1\otimes w_2))} \\
&= \ip{h_1\otimes f_1\otimes S_1}{g_2\otimes V(g_1\otimes w_1\otimes w_2)}
= \ip{h_1\otimes S_1}{g_2 \otimes w_1 \otimes \pi_2(f_1g_1)(w_2)} \\
&= \ip{T\pi_1(f_2g_2)(w_1)}{\pi_2(f_1g_1)(w_2)},
\end{align*}
and also
\begin{align*}
&\ip{f_1\otimes f_2\otimes T}{V_{13} U_{23} (g_1\otimes g_2\otimes w_1\otimes w_2)} \\
&= \ip{f_2 \otimes V^*(f_1\otimes T)}{(\chi\otimes\id\otimes\id)
  (g_1\otimes U(g_2\otimes w_1\otimes w_2))} \\
&= \ip{h_2 \otimes f_2\otimes S_2}{g_1\otimes U(g_2\otimes w_1\otimes w_2)}
= \ip{h_2 \otimes S_2}{g_1 \otimes \pi_1(f_2g_2)(w_1)\otimes w_2} \\
&= \ip{T\pi_1(f_2g_2)(w_1)}{\pi_2(f_1g_1)(w_2)},
\end{align*}
which proves equality, as required.

Finally, let $x=x_1\otimes x_2 \in E_1 \proten_\alpha E_2$ and let $\mu=\mu_1\otimes\mu_2
\in (E_1 \proten_\alpha E_2)^*$.  By Theorem~\ref{wap_rep}, if $F\in L^\infty(X)$
is defined by
\[ \ip{F}{a} = \ip{\mu}{\pi(a)(x)} \qquad(a\in L^1(X)), \]
then $F\in\wap(L^\infty(X))$.  Now, for $a\in L^1(X)$, pick $f,g\in L^2(X)$ with
$fg=a$, so we see that
\begin{align*} \ip{\mu}{\pi(a)(x)}
&= \ip{f\otimes\mu_1\otimes\mu_2}{W(g\otimes x_1\otimes x_2)}
= \ip{f\otimes\mu_1\otimes\mu_2}{UV(g\otimes x_1\otimes x_2)}.
\end{align*}

Notice that for $k\in L^2(X)$, $\lambda_1\in E_1^*$ and $\lambda_2\in E_2^*$,
\begin{align*} \ip{k\otimes\lambda_1\otimes\lambda_2}{V(g\otimes x_1\otimes x_2)}
&= \ip{\lambda_1\otimes\lambda_2}{(\id\otimes\pi_2(kg))(x_1\otimes x_2)} \\
&= \ip{\lambda_1}{x_1} \ip{\lambda_2}{\pi_2(kg)(x_2)}. \end{align*}
For $T\in (E_1\proten_\alpha E_2)^* = \mc B_{\alpha'}(E_1,E_2^*)$, as
$E\proten_\alpha F$ is reflexive, we hence must have that
\[ \ip{k\otimes T}{V(g\otimes x_1\otimes x_2)}
= \ip{T(x_1)}{\pi_2(kg)(x_2)}. \]
Define a map $\theta: L^2(X)\otimes E_2 \rightarrow L^2(X) \otimes
E_1\proten_\alpha E_2$ by $\theta(k\otimes x) = \theta(k\otimes x_1\otimes x)$
for $k\in L^2(X)$ and $x\in E_2$ on elementary tensors.  A simple calculation
shows that $\theta$ extends to a contraction $L^2(X,E_2) \rightarrow
L^2(X,E_1\proten_\alpha E_2)$.  Then, for $\tau\in L^2(X,E_2)$, we have that
\[ \ip{k\otimes\lambda_1\otimes\lambda_2}{\theta(\tau)}
= \ip{\lambda_1}{x_1} \ip{k\otimes\lambda_2}{\tau}, \]
and so, similarly,
\[ \ip{k\otimes T}{\theta(\tau)} = \ip{k\otimes T(x_1)}{\tau}. \]
Let $\hat V\in L^\infty(X)\vnten\mc B(E_2)$ be defined by $\pi_2$.
Then
\[ \ip{k\otimes T}{V(g\otimes x_1\otimes x_2)}
= \ip{k\otimes T(x_1)}{\hat V(g\otimes x_2)}
= \ip{k\otimes T}{\theta \hat V(g\otimes x_2)}. \]
Thus $V(g\otimes x_1\otimes x_2)$ is in the image of $\theta$, being
equal to $\theta \hat V(g\otimes x_2)$.

Again, we use Sweedler notation, so by the previous paragraph, we may suppose
that $V(g\otimes x_1\otimes x_2) = h \otimes x_1\otimes y_2$.  Then, for
$k\in L^2(X)$, $\lambda_1\in E_1^*$ and $\lambda_2\in E_2^*$, we see that
\begin{align*} \ip{k\otimes\lambda_1\otimes\lambda_2}{h \otimes x_1\otimes y_2}
&= \ip{\lambda_1\otimes\lambda_2}{(\id\otimes\pi_2(kg))(x_1\otimes x_2)}
= \ip{\lambda_1}{x_1} \ip{\lambda_2}{\pi_2(kg)(x_2)}.
\end{align*}
Hence we have that
\[ \ip{k\otimes\lambda_2}{h\otimes y_2}
= \ip{\lambda_2}{\pi_2(kg)(x_2)}. \]
Finally, we have that
\begin{align*} \ip{\mu}{\pi(a)(x)}
&= \ip{f\otimes\mu_1\otimes\mu_2}{U(h\otimes y_1\otimes y_2)}
= \ip{\mu_1\otimes\mu_2}{(\pi_1(fh)\otimes\id)(y_1\otimes y_2)} \\
&= \ip{\mu_2}{y_2} \ip{F_1}{fh}
= \ip{F_1 f \otimes \mu_2}{h\otimes y_2} \\
&= \ip{\mu_2}{\pi_2((F_1 f)g)(x_2)}
= \ip{F_2}{(F_1 f)g} = \ip{F_2}{F_1 fg} = \ip{F_1 F_2}{a}.
\end{align*}
So we conclude that $F_1 F_2 = F \in \wap(L^1(X))$, showing that $\wap(L^1(X))$
is an algebra.
\end{proof}

\begin{corollary}
Let $G$ be a locally compact group.  Then $\wap(M(G))$ is a C$^*$-subalgebra of
$M(G)^* = C_0(G)^{**}$.
\end{corollary}

\bigskip
\noindent\textbf{Author's address:}
\parbox[t]{5in}{School of Mathematics,\\
University of Leeds,\\
Leeds LS2 9JT\\
United Kingdom}

\smallskip
\noindent\textbf{Email:} \texttt{matt.daws@cantab.net}

\end{document}